\documentclass[a4paper, 12pt]{article}

\usepackage{amsfonts, amsmath, amsthm}
\usepackage{latexsym}

\theoremstyle{plain}
\newtheorem{thm}{Theorem}[section]

\newtheorem{lem}{Lemma}[section]

\theoremstyle{definition}
\newtheorem{dfn}{Definition}[section]

\theoremstyle{remark}
\newtheorem{rmk}{Remark}[section]

\makeatletter
 
 \@addtoreset{equation}{section}
\makeatother

\newcommand{\R}{\mathbb{R}}

\newcommand{\vomega}{\mbox{\boldmath$\omega$}}
\newcommand{\pa}{\partial}
\newcommand{\eps}{\varepsilon}

\DeclareMathOperator{\Real}{\rm Re}
\DeclareMathOperator{\Imag}{\rm Im}
\DeclareMathOperator{\supp}{\rm supp}

\begin{document}
\title{
Global small amplitude solutions for \\
 two-dimensional nonlinear Klein-Gordon systems \\
 in the presence of mass resonance
}  

\author{
         Yuichiro Kawahara
         \thanks{
             Doshisha Junior and Senior High School. 
             Iwakura, Sakyo-ku, Kyoto 606-8558, Japan. 
             (E-mail: {\tt yukawa@js.doshisha.ac.jp})
         } 
   \and  
         Hideaki Sunagawa
         \thanks{
             Department of Mathematics, Graduate School of Science, 
             Osaka University. 
             Toyonaka, Osaka 560-0043, Japan. 
             (E-mail: {\tt sunagawa@math.sci.osaka-u.ac.jp})
         }
 } 
 
\date{ }
\maketitle

\noindent{\bf Abstract:}\ 
We consider a nonlinear system of two-dimensional Klein-Gordon equations 
 with  masses $m_1$, $m_2$ satisfying the resonance relation $m_2=2m_1>0$. 
 We introduce a structural condition on the nonlinearities under which 
 the solution exists globally in time and decays at the rate $O(|t|^{-1})$ 
 as $t \to \pm \infty$ in $L^{\infty}$. 
 In particular, our new condition includes the Yukawa type interaction, 
 which has been excluded from the {\em null condition} in the sense of 
 J.-M.Delort, D.Fang and R.Xue (J.Funct.Anal.{\bf 211}(2004), 288--323).\\
 
\noindent{\bf Key Words:}\ 
Nonlinear Klein-Gordon equations; 
Mass resonance; Global solution. \\

\noindent{\bf 2000 Mathematics Subject Classification:}\ 
35L70, 35B40, 35L15\\

\section{Introduction} 
This paper is intended to be a continuation of the papers 
\cite{su1}, \cite{su2}, \cite{su3}, which are concerned with 
large time behavior of small solutions to the Cauchy problem for 
a nonlinear system of Klein-Gordon equations 
in $(t,x) \in \R^{1+2}$:
\begin{align}
 \left\{
  \begin{array}{l}
   (\Box + m_1^2) u_1 = F_1(u,\pa u), \\
   (\Box + m_2^2) u_2 = F_2(u,\pa u),
  \end{array}
 \right. 
 \label{1_1}
\end{align}
where $\Box=\pa_t^2-\pa_1^2-\pa_2^2$, $\pa=(\pa_0,\pa_1,\pa_2)$ with 
$\pa_0=\pa_t=\pa/\pa t$, $\pa_{j}=\pa/ \pa x_j$ for $j=1,2$, 
while $u=(u_j)_{j=1,2}$ is an $\R^2$-valued unknown function and 
$\pa u=(\pa_a u_j)_{{\stackrel{\scriptstyle j=1,2}{\scriptstyle a=0,1,2}}}$ 
is its first order derivative ($\R^{2\times 3}$-valued). 
The masses $m_1$, $m_2$ are supposed to be positive constants. 
Without loss of generality, we may assume that $m_1\le m_2$ throughout this 
paper. The nonlinear term $F_j =F_j(v,w)$ is a $C^{\infty}$ function of 
$(v,w) \in \R^2\times \R^{2\times 3}$ which vanishes at quadratic order at 
the origin, that is, 
$$
 F_j(v,w)=O((|v|+|w|)^2) \quad \mbox{as} \quad (v,w)\to (0,0).
$$
For simplicity, the initial data are supposed to be of the form 
\begin{align}
  u_j(0,x)=\eps f_j(x),\ \ \pa_t u_j(0,x)=\eps g_j(x), 
  \qquad  x \in \R^2,\ j=1,2 
\label{1_2}
\end{align}
with a small parameter $\eps>0$ and $C_0^{\infty}$ functions $f_j$, $g_j$.

From a perturbative point of view, quadratic nonlinear Klein-Gordon systems 
on $\R^2$ are of special interest because 
ratio of the masses and the structure of the nonlinearities play essential 
roles when one considers large time behavior of the solutions. 
Let us recall known results briefly. 
In the case of $m_2\ne 2m_1$ (which will be referred to as the 
{\em non-resonant} case), it is shown in \cite{su1}, \cite{tsu} that 
the solution $u(t)$ for (\ref{1_1})--(\ref{1_2}) exists globally without any 
structural restrictions of $F_1$, $F_2$ if $\eps$ is sufficiently small. 
Moreover, $u(t)$ is asymptotically free (in the sense that we can find a 
solution $u^{\pm}(t)$ of the homogeneous linear Klein-Gordon equations 
such that $u(t)$ tends to $u^{\pm}(t)$ as $t \to \pm \infty$ in the energy 
norm) and 
satisfies the following time decay estimate for all $p\in [2,\infty]$:  
\begin{align}
 \sum_{|I|\le 1} \|\pa_{t,x}^{I} u(t,\cdot)\|_{L^p(\R_x^2)}
 \le 
 C\eps (1+|t|)^{-\left(1-\frac{2}{p}\right)}
  \qquad (t \in \R)
   \label{1_3}
\end{align}
with some positive constant $C$ which is independent of $\eps$. 
Remember that this decay rate is same as that for the linear case. 
On the other hand, the above assertion fails to hold in the {\em resonant} 
case (i.e., the case where $m_2 = 2m_1$) 
because of counterexamples due to \cite{su2}, \cite{su3}, \cite{taf} etc. 
One of the  simplest example is
\begin{align*}
\left\{\begin{array}{l}
F_1=0 \\
F_2=u_1^2.
\end{array}\right.
\end{align*}
For this nonlinearity, we can choose $f_j$, $g_j \in C_0^{\infty}(\R^2)$ 
and positive constants $C$, $T$ such that the solution $u(t)$ for 
(\ref{1_1})--(\ref{1_2}) satisfies 
 $$
  \sum_{|I|\le 1} \|\pa_{t,x}^{I} u(t,\cdot)\|_{L^2} \ge C \eps^2 \log |t| 
  \qquad (|t| \ge T)
 $$
however small $\eps>0$  is, whence the estimate (\ref{1_3}) is violated. 
Thus we need to put some structural condition on the nonlinearities in order 
to obtain global solutions for (\ref{1_1})--(\ref{1_2}) satisfying (\ref{1_3}) 
in the resonant case. This is what we are going to address here. 
A sufficient condition on the nonlinearities is introduced by 
Delort--Fang--Xue \cite{dfx}, called the {\em null condition}, 
which admits a global solution for 
(\ref{1_1})--(\ref{1_2}) in the resonant case. They also give an asymptotic 
profile of the solution, from which the decay estimate  (\ref{1_3}) follows 
immediately. However, their condition is not optimal since it does not cover 
some important cases. For instance, 
\begin{align}
\left\{\begin{array}{l}
F_1=u_1u_2 \\
F_2=u_1^2
\end{array}\right.
\label{1_4}
\end{align}
is excluded from their condition, while the system (\ref{1_1}) with 
the nonlinearity (\ref{1_4}) can be viewed as a simplified model for some 
physical systems, such as Dirac-Klein-Gordon system, Maxwell-Higgs system, 
and so on. Someone may call this type of interaction the {\em Yukawa type} 
one (see e.g., \cite{hnb}, \cite{tsu} and the references therein). 

Our aim in this paper is to give a new sufficient condition on the 
nonlinearities which includes (\ref{1_4}). 
Under this condition, we will show that the solution for 
(\ref{1_1})--(\ref{1_2}) exists globally in time and it enjoys time decay 
property (\ref{1_3}) even in the resonant case.

\section{Main result} 

In order to state the result, let us introduce several notations. 
For $j=1,2$, denote by $Q_j$ the quadratic homogeneous part of the nonlinear 
term $F_j$, that is, 
$$
 Q_j(v,w)=\lim_{\lambda \downarrow 0} \lambda^{-2} F_j(\lambda v, \lambda w) 
$$
for $(v,w) \in \R^2\times \R^{2\times 3}$. 
Roughly saying, $Q_j(u,\pa u)$ gives the main part of the 
nonlinearity while $F_j(u,\pa u)-Q_j(u,\pa u)$ is regarded as a cubic or 
higher order remainder if we are interested in small amplitude solutions. 
Next we set 
\begin{align*}
 \mathbb{H} =\{\vomega=(\omega_0,\omega_1,\omega_2) \in \R^3\, :\, 
  \omega_0^2-\omega_1^2-\omega_2^2=1  \}
\end{align*}
and 
\begin{align}
 \Phi_j(\vomega)
 =\int_{0}^{1} 
  Q_j \bigl( V(\theta), W(\vomega,\theta) \bigr)\, 
 e^{-2 \pi i j\theta}\, d\theta
 \label{2_1}
\end{align}
for $\vomega \in \mathbb{H}$, where 
$V(\theta) = \bigl(\cos 2 \pi k\theta \bigr)_{k=1,2}$, 
$W(\vomega,\theta) = 
  \bigl( -\omega_a m_k \sin 2\pi k\theta 
  \bigr)_{{\stackrel{\scriptstyle k=1,2}{\scriptstyle a=0,1,2}}}$ 
and $i=\sqrt{-1}$. 
Note that $\Phi_j$ can be explicitly computed only from $m_1$, $m_2$ and 
$F_j$. With these $\Phi_1(\vomega)$ and $\Phi_2(\vomega)$, 
we introduce the following two conditions: 
\begin{quote}
(a)\ 
Both $\Phi_1(\vomega)$ and $\Phi_2(\vomega)$ vanish identically on 
$\mathbb{H}$. 
\end{quote}
\begin{quote}
(b)\ 
The real part of the product $\Phi_1(\vomega) \Phi_2(\vomega)$ is uniformly 
positive on $\mathbb{H}$, while the imaginary part of 
$\Phi_1(\vomega) \Phi_2(\vomega)$ vanishes identically on $\mathbb{H}$. 
\end{quote}
Our  main result is the following theorem.
\begin{thm}\label{thm2_1}
Let $m_2=2m_1>0$. 
 Suppose that either the condition (a) or (b) is satisfied. 
Then (\ref{1_1})--(\ref{1_2}) admits a unique global classical solution for 
sufficiently small $\eps$. Moreover, for all $p \in [2,\infty]$, 
the solution $u(t)$ satisfies (1.3), i.e., 
\begin{align*}
 \sum_{|I|\le 1} \|\pa_{t,x}^{I} u(t,\cdot)\|_{L^p(\R_x^2)}
 \le 
 C\eps (1+|t|)^{-\left(1-\frac{2}{p}\right)}
\end{align*}
with some positive constant $C$ which does not depend on $\eps$. 
\end{thm}

\begin{rmk}
The condition (a) is equivalent to the null condition in the sense of 
\cite{dfx}. On the other hand, the condition (b) is completely 
new, as far as the authors know. (\ref{1_4}) is a typical example 
of the nonlinearity which is excluded from (a) but included in (b). 
As our proof below suggests, it may be reasonable to conjecture that 
the solution may {\em not} be asymptotically free under the condition (b)  
(while it is possible to prove that the solution is asymptotically free under 
the condition (a); see \cite{kos} for the detail). 
This problem will be discussed in a future work.
\end{rmk}

\begin{rmk}
Our main result remains valid for quasilinear systems if 
the definition of $\Phi_j$ is slightly modified and 
a suitable hyperbolicity assumption is imposed on $F_j$.\\
\end{rmk}

The rest of this paper is organized as follows. In Section 3 we make some 
reduction of the problem along the idea of \cite{de}, \cite{dfx} with a 
slight modification. 
Section 4 is devoted to the derivation of some energy inequalities. 
In Section 5 we specify the worst contribution of the nonlinearities in the 
resonant case. Section 6 describes a lemma on some ordinary differential 
equations, which reveals the role of our condition imposed on $\Phi_j$. 
After that, we get an a priori estimate in Section 7, 
from which global existence follows immediately. 
Finally, in Section 8, the time decay estimate (\ref{1_3}) is derived. 
In what follows, several positive constants appearing in estimates will 
be denoted by the same letter C, which may vary from line to line.

\section{Reduction of the problem} 
In the following, we restrict ourselves to the forward Cauchy problem ($t>0$)
since the backward problem can be treated in the same way. 
Also, we shall neglect the higher order terms of $F_j$ (i.e. we assume 
$F_j=Q_j$) to make the essential idea clearer.

Let $K$ be a positive constant which satisfies
$$
 \supp f_j \cup \supp g_j \subset \bigl\{ x\in \R\; :\; |x|\leq K \bigr\}
$$
and let $\tau_0$ be a fixed positive number strictly greater than $1+2K$. 
We start with the fact that we may treat the problem as if the Cauchy data are 
given on the upper branch of the hyperbola 
$$
 \bigl\{ (t,x) \in \R^{1+2}\; :\;  (t+2K)^2-|x|^2 = \tau_0^2,\ t>0 \bigr\}
$$ 
and it is sufficiently smooth, small, compactly-supported. 
This is a consequence of the classical local existence theorem and the 
finite speed of propagation 
(see e.g., \cite[Proposition 1.4]{de} or \cite[Proposition 1.1.4]{dfx} for 
the detail). 
Next, let us introduce the hyperbolic coordinate 
$(\tau, z) \in [\tau_0,\infty)\times \R^2$ in the interior of the light cone, 
i.e., 
$$
 t+2K = \tau \cosh |z|,\ \ 
  x_1 = \tau \frac{z_1 }{|z|}\sinh |z|, \ \ 
  x_2 = \tau  \frac{z_2}{|z|} \sinh |z|
$$ 
for $|x| < t+2K$. 
Then, with the auxiliary expression 
$z_1=\rho\cos \theta$, $z_2=\rho \sin \theta$, 
we see that 
$$
 \begin{pmatrix}
 \pa_{0} \\ \pa_{1} \\ \pa_{2}
 \end{pmatrix}
 =
 \begin{pmatrix}
 \pa_{t} \\ \pa_{x_1} \\ \pa_{x_2}
 \end{pmatrix}
 =
 \begin{pmatrix}
 \cosh \rho & -\sinh \rho & 0\\
 -\sinh \rho \cos \theta & \cosh \rho \cos \theta & -\sin \theta \\
 -\sinh \rho \sin \theta & \cosh \rho \sin \theta & \cos \theta
 \end{pmatrix}
 \begin{pmatrix}
\pa_{\tau} \\ \frac{1}{\tau}\pa_{\rho} \\ \frac{1}{\tau \sinh \rho}\pa_{\theta}
 \end{pmatrix}, 
$$
whence
\begin{align}  
\pa_a= 
  \omega_a(z)\pa_{\tau}  + \frac{1}{\tau}\sum_{j=1}^{2} \eta_{aj}(z)\pa_{z_j}
 \label{3_1}
\end{align}
for $a=0,1,2$, where
\begin{align*} 
 \vomega(z)=
  \begin{pmatrix}
   \omega_0(z) \\ \omega_1(z) \\ \omega_2(z)
  \end{pmatrix}
  =
  \begin{pmatrix}
   \cosh\rho \\ -\sinh\rho \cos \theta\\ -\sinh\rho \sin \theta
  \end{pmatrix},
\end{align*} 
\begin{align*} 
\begin{pmatrix}
 \eta_{01}(z) & \eta_{02}(z)\\
 \eta_{11}(z) & \eta_{12}(z)\\
 \eta_{21}(z) & \eta_{22}(z)\\
\end{pmatrix}
 =\begin{pmatrix}
- \sinh \rho \cos \theta
&
- \sinh \rho \sin \theta
\\
 \cosh\rho \cos^2 \theta+ \frac{\rho}{\sinh \rho} \sin^2 \theta
&
\left(\cosh\rho - \frac{\rho}{\sinh \rho} \right) \cos \theta \sin \theta
\\
\left(\cosh\rho - \frac{\rho}{\sinh \rho} \right) \cos \theta \sin \theta
&
\cosh\rho \sin^2 \theta+ \frac{\rho}{\sinh \rho}\cos^2 \theta
\end{pmatrix}.
\end{align*} 
Remark that $\omega_a(z)$ and $\eta_{bj}(z)$ can be regarded as 
$C^{\infty}$ functions of $z \in \R^2$ which satisfy 
$$
 |\omega_a(z)| + |\eta_{bj}(z)| \le C e^{|z|} 
$$
for $a,b =0,1,2$ and $j=1,2$. 
Moreover, $\vomega(z) \in \mathbb{H}$ for all $z \in \R^2$. 
Also we observe that 
\begin{align*}
 \Box u 
  =  \frac{1}{\tau} 
     \left(
      \pa_{\tau}^2 -\frac{1}{\tau^2}\Lambda_0
     \right) 
     \bigl( \tau u \bigr),
\end{align*}
where 
\begin{align}
\Lambda_0 =\pa_{\rho}^2+\frac{\cosh \rho}{\sinh \rho}\pa_{\rho} 
            +
            \frac{1}{\sinh^2 \rho}\pa_{\theta}^2.
\label{3_2}            
\end{align}
Next we introduce a weight function $\chi(z) =e^{-\kappa\langle z \rangle}$ 
with a large parameter $\kappa$, 
where $\langle z \rangle=\sqrt{1+|z|^2}$.  
(In fact, we shall not always need the explicit form of $\chi$, but only the 
properties that $\chi$ is smooth, radial, as well as the estimates 
$0 < \chi(z) \leq C_0 e^{-\kappa |z|}$ and 
$|\pa_z^{I} \chi (z)| \leq C_{I} \chi(z)$ for any multi-indices $I$ 
with some constants $C_{I}$. Another choice for such $\chi(z)$ may be 
$\frac{1}{\cosh(\kappa|z|)}$, as was done by Delort et al. in \cite{de}, 
\cite{dfx}. We also note that $\kappa \ge 6$ is enough for our 
purpose.)  With this weight function, 
let us define the new unknown function $v_j(\tau,z)$ by 
\begin{align*}
 u_j(t,x) = \frac{\chi(z)}{\tau} v_j(\tau,z).
\end{align*}
Then we see that $v=(v_1,v_2)$ satisfies
\begin{align*}
 \left( \pa_{\tau}^2 - \frac{1}{\tau^2}\Lambda + m_j^2 \right)v_j 
 = 
 \tilde{Q_j}(\tau,z,v, \pa_{\tau,z} v)
\end{align*}
if $u=(u_1,u_2)$ solves (1.1), where $\Lambda$ is defined by 
\begin{align*}
 \Lambda v  = e^{\kappa \langle z \rangle}  
        \Lambda_0 \bigl( e^{-\kappa \langle z \rangle} v \bigr)
\end{align*}
and 
\begin{align} 
 \tilde{Q}_j(\tau,z,v, \pa_{\tau,z} v)= 
 \frac{\chi(z)}{\tau}Q_j(v,\vomega(z)\pa_{\tau}v)
 +
 \sum_{\nu=0}^{1}
 \sum_{
  {\stackrel{\scriptstyle 1\le k,l \le 2}{\scriptstyle |I|\le 1, |J|\le \nu}}
 }
 \frac{q_{\nu jklIJ}(z)}{\tau^{2+\nu}} 
 \pa_{z}^I v_k \cdot \pa_{\tau}^{1-\nu}\pa_{z}^J v_l
 \label{3_3}
\end{align}
with some $q_{\nu jklIJ} \in C^{\infty}(\R^2)$ satisfying
$$
 |\pa_z^L q_{\nu jklIJ}(z)|\le C_L e^{(2-\kappa)|z|}
$$
for any multi-index $L$.

At last, the original problem (\ref{1_1})--(\ref{1_2}) is reduced to 
\begin{align}
 \left\{ 
 \begin{array}{cc}
  \left( \pa_{\tau}^2 -\frac{1}{\tau^2}\Lambda +m_j^2 \right)v_j
  = 
  \tilde{Q}_j(\tau,z,v, \pa_{\tau,z} v), 
  & \tau > \tau_0,\ z \in \R^2,\\ 
  (v_j, \pa_{\tau} v_j)|_{\tau=\tau_0}=(\eps \tilde{f}_j, \eps \tilde{g}_j)
  & z \in \R^2,
 \end{array}
 \right.
 \label{3_4}
\end{align}
where $\tilde{f}_j$ and $\tilde{g}_j$ are $C^{\infty}$
functions of $z \in \R^2$ with compact support. \\

\section{Commuting vector fields and energy inequalities} 
In this section,  we will derive a kind of energy inequalities for the 
operator 
$$
 P_m =\pa_{\tau}^2 -\frac{1}{\tau^2}\Lambda +m^2
$$
with $m>0$ which will be needed in Section 7. 
For this purpose it is helpful to introduce the following function class.
\begin{dfn}
Let $\nu \in \R$. We denote by $\mathcal{S}^{\nu}$ the space of 
$C^{\infty}$ functions $a(z)$ defined on $\R^2$ satisfying
$$
\sup_{z\in \R^2} 
\left(\frac{|\pa_z^{I} a(z)|}{\langle z \rangle^{\nu-|I|}} \right)
<\infty 
$$
for any multi-index $I$.
\end{dfn}

We start with splitting $\Lambda$ into three parts: 
$\Lambda = \Lambda_0 +\Lambda_1 +\Lambda_2$, 
where $\Lambda_0$ is defined by (\ref{3_2}) and 
\begin{align*}
 &\Lambda_1 
 = 
 -2\kappa \frac{\rho}{\langle \rho \rangle} \pa_{\rho}
 =
 -2\kappa \sum_{j=1}^{2}\frac{z_j}{\langle z \rangle} \pa_{z_j},\\
\quad 
&\Lambda_2 =
 \frac{\kappa^2 |z|^2\langle z \rangle-\kappa}{\langle z \rangle^3}
 +
 \frac{\kappa |z| \cosh |z|}{\langle z \rangle \sinh |z|}.
\end{align*}
Note that we can rewrite $\Lambda_0$ as 
\begin{align*} 
 \Lambda_0 v =\frac{1}{\sqrt{\mathcal{G}(z)}}
            \sum_{i,j=1}^{2}\pa_{z_i} 
            \left(\sqrt{\mathcal{G}(z)}g^{ij}(z) \pa_{z_j} v \right)
\end{align*} 
when we put 
\begin{align*} 
 \begin{pmatrix} g^{11}(z)& g^{12}(z)\\ g^{21}(z) & g^{22}(z) \end{pmatrix}
 = 
 \begin{pmatrix} 
 1 &0 \\
 0 & 1
 \end{pmatrix}
 - 
 \left(\frac{1}{|z|^2} - \frac{1}{\sinh^2 |z|}\right)
 \begin{pmatrix} 
 z_2^2 & -z_1z_2 \\
 -z_1z_2 & z_1^2
 \end{pmatrix}
\end{align*} 
and
\begin{align*} 
 \mathcal{G}(z)
  = \left(\frac{\sinh |z|}{|z|} \right)^2.
\end{align*} 
We observe that
\begin{align} 
 \sum_{j,k=1}^{2}g^{jk}(z)\zeta_j\zeta_k 
 = 
 \left|\frac{z}{|z|}\cdot \zeta\right|^2
 +
 \frac{1}{\sinh^2|z|}\left|z \wedge \zeta\right|^2 \ge 0
 \label{4_1}
\end{align} 
for $\zeta=(\zeta_1,\zeta_2) \in \R^2$, and that 
\begin{align*} 
 \mathcal{G}(z)
  =\det(g^{jk}(z))_{1\le j,k \le 2}^{-1}.
\end{align*} 
Next, let us introduce the  vector fields 
\begin{align*} 
 &\Gamma_1
 =(t+2K)\pa_{x_1}+x_1\pa_t
 =(\cos \theta) \pa_{\rho} -\frac{\sin \theta}{\tanh \rho} \pa_{\theta}, 
\\
 &\Gamma_2
 =(t+2K)\pa_{x_2}+x_2\pa_t
 =(\sin \theta) \pa_{\rho} +\frac{\cos \theta}{\tanh \rho} \pa_{\theta}, 
\\
 &\Gamma_3 
 = -x_2\pa_{x_1}+x_1\pa_{x_2}
 =\pa_{\theta}.
\end{align*} 
In what follows, we write $|I|=I_1+I_2+I_3$ and 
$\Gamma^I=\Gamma_1^{I_1}\Gamma_2^{I_2}\Gamma_3^{I_3}$ 
for a multi-index $I=(I_1,I_2,I_3)$. 
We can immediately check that
$$
  [\Gamma_1,\Gamma_2]=\Gamma_3,\ \ 
  [\Gamma_1, \Gamma_3]=\Gamma_2, \ \ 
  [\Gamma_2, \Gamma_3]=\Gamma_1,
$$
where $[\cdot,\cdot]$ denotes the commutator. 
Another important thing is that $\Gamma_1$, $\Gamma_2$  are written as linear 
combinations of $\pa_{z_1}$, $\pa_{z_2}$ with $\mathcal{S}^1$-coefficients, 
while $\pa_{z_1}$, $\pa_{z_2}$ are written as linear 
combinations of $\Gamma_1$, $\Gamma_2$ with $\mathcal{S}^0$-coefficients. 
More precisely, we have
\begin{align*}
 \Gamma_j=\sum_{k=1}^{2} c_{jk}(z) \pa_{z_k}
\end{align*}
for $j=1,2$ and 
\begin{align*}
 \pa_{z_k}=\sum_{l=1}^{2} \tilde{c}_{kl}(z) \Gamma_l
\end{align*}
for $k=1,2$, where
\begin{align*}
 &\begin{pmatrix} c_{11}(z)& c_{12}(z)\\ c_{21}(z) & c_{22}(z) \end{pmatrix}
 = 
 \begin{pmatrix} 
  \cos \theta & -\sin \theta\\ \sin \theta & \cos \theta
 \end{pmatrix}
 \begin{pmatrix} 
  1 & 0\\ 0 & \frac{\rho}{\tanh \rho}
 \end{pmatrix}
 \begin{pmatrix} 
  \cos \theta & \sin \theta\\ -\sin \theta & \cos \theta
 \end{pmatrix},\\
 &\begin{pmatrix} 
  \tilde{c}_{11}(z)& \tilde{c}_{12}(z)\\ \tilde{c}_{21}(z) & \tilde{c}_{22}(z) 
 \end{pmatrix}
 = 
  \begin{pmatrix} 
  \cos \theta & -\sin \theta\\ \sin \theta & \cos \theta
 \end{pmatrix}
 \begin{pmatrix} 
  1 & 0\\ 0 & \frac{\tanh \rho}{\rho}
 \end{pmatrix}
 \begin{pmatrix} 
  \cos \theta & \sin \theta\\ -\sin \theta & \cos \theta
 \end{pmatrix}.
\end{align*}

As for the commutation relation 
between $P_m$ and $\Gamma_j$'s, we have the following: 
\begin{lem} \label{lem4_1}
For any multi-index $I$, we have
$$
 [P_m,\Gamma^I]
 =
 \frac{1}{\tau^2}\sum_{|J| \le |I|} h_{IJ}(z)\Gamma^J
$$ 
with some  $h_{IJ}\in \mathcal{S}^{0}$.
\end{lem}

Proof: 
First we note  that 
$$
 [\Box +m^2, \Gamma^I]=0
$$
for any multi-index $I$, and that 
\begin{align*}
\Box+m^2
=
\pa_{\tau}^2+\frac{2}{\tau}\pa_{\tau} -\frac{1}{\tau^2}\Lambda_0+m^2
=
P_m +\frac{1}{\tau^2}(\Lambda_1+\Lambda_2). 
\end{align*}
So we have 
\begin{align*}
 [P_m,\Gamma^I]
 &= 
 -\frac{1}{\tau^2}[\Lambda_1+\Lambda_2,\Gamma^I]\\
 &= 
 \frac{1}{\tau^2}
 \sum_{j=1}^{2}[p_j\Gamma_j,\Gamma^I]
 -\frac{1}{\tau^2}[\Lambda_2,\Gamma^I],
\end{align*}
where $p_j(z)=2\kappa\frac{z_j}{\langle z\rangle} \in \mathcal{S}^0$. 
By induction on $I$, we have the desired conclusion.
\qed\\

Now, we turn to the energy inequalities for the operator $P_m$ 
which we need. For $s\in \mathbb{Z}_{\ge 0}$, we introduce the energy $E_s$ 
as follows:
\begin{align*}
 E_s(\tau;v,m) 
 = \sum_{|I|\le s} \frac{1}{2}
 \int_{\R^2} \left(
 \bigl(\pa_{\tau} \Gamma^I v \bigr)^2 
 + \frac{1}{\tau^2}\sum_{j,k=1}^{2}g^{jk}(z) 
   \bigl(\pa_{z_j}\Gamma^I v \bigr) \bigl(\pa_{z_k}\Gamma^I v \bigr)
 + m^2 \bigl( \Gamma^{I} v \bigr)^2 
 \right) \sqrt{\mathcal{G}(z)} d z.
\end{align*}
We also introduce the norm $\| \cdot \|_{(s)}$ by
$$
 \|v\|_{(s)}:=
 \sum_{|I|\le s} \|\Gamma^I v\|_{L^2(\R^2;\sqrt{\mathcal{G}(z)}dz)}.
$$

\begin{lem} \label{lem4_2}
For $s \in \mathbb{Z}_{\ge 0}$, we have 
\begin{eqnarray}
 \frac{d}{d \tau} E_s(\tau;v,m) 
 \leq 
 \left(\frac{2\kappa}{\tau} + \frac{C}{\tau^2}\right) E_s(\tau;v,m)
 + CE_s(\tau;v,m)^{1/2} \|P_m v(\tau) \|_{(s)}
 \label{4_2}
\end{eqnarray}
and 
\begin{eqnarray}
 \frac{d}{d \tau} E_s(\tau;v,m) 
 \leq 
 \frac{C}{\tau^{2}} E_{s+1}(\tau;v,m) 
 + CE_s(\tau;v,m)^{1/2} \|P_m v(\tau) \|_{(s)}. 
\label{4_3}
\end{eqnarray}

\end{lem}

Proof: 
First we consider the case of $s=0$. As usual, we compute 
\begin{align*}
 &\frac{d}{d\tau}E_0(\tau;v,m)\\
 &=
 \int_{\R^2}
  \biggl(
  (\pa_{\tau}v) \pa_{\tau}^2v
  +
  m^2v \pa_{\tau}v
  + 
  \frac{1}{\tau^2}\sum_{j,k=1}^{2}g^{jk}(z)(\pa_{z_k}v)\pa_{\tau}\pa_{z_j}v
  -\frac{2}{\tau^3}
  \sum_{j,k=1}^{2}g^{jk}(z)(\pa_{z_j}v)\pa_{z_k}v
 \biggr)\sqrt{\mathcal{G}(z)}dz\\
 &\le
 \int_{\R^2}
  (\pa_{\tau}^2v  +  m^2v) (\pa_{\tau}v) \sqrt{\mathcal{G}(z)}
  - 
  \frac{1}{\tau^2}\sum_{j,k=1}^{2}
  \pa_{z_j}\Bigl( \sqrt{\mathcal{G}(z)}g^{jk}(z)\pa_{z_k}v \Bigr)
  (\pa_{\tau}v)
  dz\\
&=
 \int_{\R^2}
  \Bigl(P_m v-\frac{1}{\tau^2} \Lambda_1v -\frac{1}{\tau^2} \Lambda_2v \Bigr)
  (\pa_{\tau}v) \sqrt{\mathcal{G}(z)}dz\\
&\le
 \|P_m v\|_{(0)} E_0(\tau;v,m)^{1/2}
 +
 \frac{1}{\tau^{1+l}} \int_{\R^2}
 \frac{|\Lambda_1 v|}{\tau^{1-l}} |\pa_{\tau} v|\sqrt{\mathcal{G}(z)}dz
 +
 \frac{C}{\tau^2} E_0(\tau;v,m)
\end{align*}
for $l=0,1$. 
We shall estimate the second term differently according to $l=0$ or $l=1$. 
In the case of $l=0$, from the relations 
$$
 \left| \Lambda_1 v \right| 
 =
 \left| 2\kappa \frac{ |z|}{\langle z \rangle} \pa_{\rho} v \right|
 \le 
 2\kappa \left| \pa_{\rho} v \right|
$$
and
$$
 \sum_{j,k=1}^{2}g^{jk}(z)(\pa_{z_j}v)(\pa_{z_k}v)
 =
 \left| \pa_{\rho} v \right|^2 
 + 
 \frac{1}{\sinh^2|z|}\left| \pa_{\theta} v \right|^2
 \ge 
 \left| \pa_{\rho} v \right|^2
$$
it follows that 
\begin{align*}
 \frac{1}{\tau}
 \int_{\R^2}\frac{|\Lambda_1v|}{\tau}|\pa_{\tau} v|\sqrt{\mathcal{G}(z)}dz
 &\le 
 \frac{2\kappa}{\tau} 
 \int_{\R^2} \frac{|\pa_{\rho}  v|}{\tau}|\pa_{\tau} v| \sqrt{\mathcal{G}(z)}dz
 \\
 &\le 
 \frac{\kappa}{\tau} \int_{\R^2}
  \left(\frac{|\pa_{\rho} v|^2}{\tau^2}+ |\pa_\tau v|^2 \right)
 \sqrt{\mathcal{G}(z)}dz\\
 &\le 
 \frac{\kappa}{\tau}  \int_{\R^2}\left(
 \frac{1}{\tau^2}\sum_{j,k=1}^{2}g^{jk}(z)(\pa_{z_j}v)(\pa_{z_k}v)
 +
  |\pa_\tau v|^2 
 \right)\sqrt{\mathcal{G}(z)}dz\\
 &\le 
 \frac{2\kappa}{\tau} E_0(\tau;v,m),
\end{align*}
which gives us $(\ref{4_2})$ with $s=0$. On the other hand, using the relation 
$$
 |\Lambda_1 v||\pa_{\tau} v| 
 =
 \left|2\kappa \sum_{j=1}^{2}\frac{z_j}{\langle z \rangle }\Gamma_j v\right|
 |\pa_{\tau}v|
 \le 
 \frac{\kappa}{m}\bigl( m^2|\Gamma v|^2+|\pa_{\tau}v|^2 \bigr),
$$
we have 
\begin{align*}
 \frac{1}{\tau^2}
 \int_{\R^2}|\Lambda_1 v||\pa_{\tau} v|\sqrt{\mathcal{G}(z)}dz
 &\le 
 \frac{C}{\tau^2} \int_{\R^2}
    \left( m^2|\Gamma v|^2+ |\pa_\tau v|^2 \right)
  \sqrt{\mathcal{G}(z)}dz\\
 &\le 
 \frac{C}{\tau^2} E_1(\tau;v,m),
\end{align*}
which yields $(\ref{4_3})$ with $s=0$. 
Next we consider the case of $s \ge 1$. 
It follows from Lemma \ref{lem4_1} that 
$$
 \sum_{|I|\le s}\|[P_m,\Gamma^I]v\|_{(0)} 
 \le \frac{C}{\tau^2} \|v\|_{(s)}
 \le \frac{C}{\tau^2} E_s(\tau;v,m)^{1/2}.
$$
Therefore
\begin{align*}
 &\frac{d}{d\tau}E_s(\tau;v,m)
 =
 \sum_{|I|\le s} \frac{d}{d\tau}E_0(\tau; \Gamma^I v, m)\\
 &\le 
  \sum_{|I|\le s} 
  \biggl\{
   \left(\frac{2\kappa}{\tau} + \frac{C}{\tau^2}\right) E_0(\tau;\Gamma^Iv,m)
   + 
   CE_0(\tau;\Gamma^I v,m)^{1/2} \|P_m \Gamma^Iv(\tau) \|_{(0)}
 \biggr\}\\
 &\le 
  \left(\frac{2\kappa}{\tau} + \frac{C}{\tau^2}\right) E_s(\tau;v,m)
   + 
   CE_s(\tau;v,m)^{1/2}
    \sum_{|I|\le s}\Bigl(
      \|\Gamma^IP_m v(\tau) \|_{(0)} +\|[P_m,\Gamma^I]v\|_{(0)}
   \Bigr)\\
 &\le
   \left(\frac{2\kappa}{\tau} + \frac{C}{\tau^2}\right) E_s(\tau;v,m)
   + 
   CE_s(\tau;v,m)^{1/2}\|P_m v(\tau) \|_{(s)}.
\end{align*}
This completes the proof of $(\ref{4_2})$. 
In the same way $(\ref{4_3})$ can be derived.
\qed

We close this section with the following lemma, which will be used in 
Section 7 to estimate quadratic terms.

\begin{lem} \label{lem4_3}
For $\kappa>9/2$ and $s \in \mathbb{Z}_{\ge 0}$, we have 
\begin{align*}
\|e^{-\kappa\langle z \rangle}\varphi \psi\|_{(s)}
 \le 
C \Bigl(
\|e^{-2|z|} \varphi\|_{L^{\infty}} \| \psi\|_{(s)}
+ 
\| \varphi \|_{(s)} \|e^{-2|z|} \psi\|_{L^{\infty}}
\Bigr),
\end{align*}
provided that the right hand side is finite.
\end{lem}

Proof: First we note that
$$
 \sum_{|I|\le s} |\Gamma^I \phi(z)|^2 \sqrt{\mathcal{G}(z)}
 \le 
 C \sum_{|I|+j\le s}  
 \bigl|
  e^{(1/2)\langle z \rangle} 
  \langle z \rangle^{|I|}\pa_z^I\pa_{\theta}^j \phi(z)
 \bigr|^2,
$$
whence
$$
 \|\phi\|_{(s)}
 \le 
 C 
 \sum_{|I|+j \le s} \bigl\| 
    e^{(1/2+\delta)\langle z \rangle} \pa_z^{I}\pa_{\theta}^{j} \phi 
    \bigr\|_{L^2(\R^2;dz)}
$$
for any $\delta>0$. By taking $\delta=\kappa -9/2$ (so that 
$1/2+\delta=\kappa-4$), we have 
\begin{align*}
 &\| e^{-\kappa\langle z \rangle} \varphi \psi \|_{(s)}\\
 &\le 
  C \sum_{|I|+j \le s} \bigl\| 
     e^{(1/2+\delta)\langle z \rangle} \pa_z^{I}\pa_{\theta}^{j}
     (e^{-\kappa\langle z \rangle} \varphi \psi) 
    \bigr\|_{L^2}\\
 &= 
  C \sum_{|I|+j \le s} \Bigl\| e^{(\kappa-4)\langle z \rangle} 
    \pa_z^{I}\bigl\{ e^{-(\kappa-4) \langle z \rangle}
    \pa_{\theta}^{j}
    (e^{-2\langle z \rangle}\varphi \cdot e^{-2\langle z \rangle} \psi) 
    \bigr\} \Bigr\|_{L^2}\\
 &\le 
  C \sum_{|I|+j \le s} 
 \|\pa_z^{I}\pa_{\theta}^{j}(
  e^{-2\langle z \rangle} \varphi\cdot 
  e^{-2\langle z \rangle} \psi)
 \|_{L^2}\\
 &\le 
  C\Bigl\{
  \|e^{-2\langle z \rangle} \varphi \|_{L^{\infty}}
  \sum_{|I|+j \le s} \| 
    \pa_{z}^{I} \pa_{\theta}^j (e^{-2\langle z \rangle} \psi)
  \|_{L^2}
  + 
   \|e^{-2\langle z \rangle} \psi\|_{L^{\infty}}
  \sum_{|I|+j \le s} \| 
   \pa_{z}^{I} \pa_{\theta}^j (e^{-2\langle z \rangle} \varphi)
  \|_{L^2}
 \Bigr\}\\
 &\le
  C\Bigl\{
   \| e^{-2|z|} \varphi \|_{L^{\infty}}
   \sum_{|I| \le s} \| \Gamma^{I} \psi  \|_{L^2}
  +  
   \|e^{-2|z|} \psi\|_{L^{\infty}}
   \sum_{|I| \le s} \| \Gamma^{I} \varphi  \|_{L^2}
  \Bigr\}\\ 
 &\le
  C\Bigl\{
  \|e^{-2|z|} \varphi\|_{L^{\infty}} \|\psi \|_{(s)}  
  + 
  \|e^{-2|z|} \psi\|_{L^{\infty}} \|\varphi\|_{(s)}
  \Bigr\}.
\end{align*}

\qed

\section{The leading part of the nonlinearity} 

The objective of this section is to extract the leading part of 
${Q}_j(v,\vomega \pa_{\tau}v)$ under some assumptions on $v$. 
What we are going to prove is the following:

\begin{lem} \label{lem5_1}
Let $m_2= 2 m_1 >0$, $\vomega =(\omega_a)_{a=0,1,2}\in \mathbb{H}$, 
$T> \tau_0>0$ and $\eps>0$. 
Suppose that $v=(v_1,v_2)$ is an $\R^2$-valued function of 
$(\tau, z) \in [\tau_0,T)\times \R^2$ which satisfies
$$
 |v_j(\tau,z)| + |\pa_{\tau}v_j(\tau,z)| \le C\eps^{1/2} e^{2|z|}, 
 \quad 
 |(\pa_{\tau}^2+m_j^2)v_j(\tau,z)|
 \le 
 \frac{C\eps^{1/2} e^{2|z|}}{ \tau}
$$
for $(\tau,z) \in  [\tau_0,T)\times \R^2$, $j=1,2$. 
Then  we have 
\begin{align}
 & 
 \left|
  \frac{e^{-i m_1 \tau} }{\tau}Q_1(v,\vomega \pa_{\tau}v)
  -
 \Biggl( 
  \frac{\Phi_1(\vomega)}{\tau} \overline{\alpha_1} \alpha_2
  +
  \pa_{\tau}\gamma_1
 \Biggr) \right|
 \le  
 \frac{C \eps \langle \vomega \rangle^2 e^{4|z|}}{\tau^2},
 \label{5_1}
 \\
 &\left|
  \frac{e^{-im_2 \tau} }{\tau}Q_2(v,\vomega \pa_{\tau}v)
 -  
 \Biggl( 
  \frac{ \Phi_2(\vomega)}{\tau}\alpha_1^2
  +
  \pa_{\tau}\gamma_2 
 \Biggr) \right|
 \le  
 \frac{C \eps \langle \vomega \rangle^2e^{4|z|}}{\tau^2},
 \label{5_2}
\end{align}
where 
$\Phi_j(\vomega)$ is given by (\ref{2_1}), $\alpha_j$ is defined by 
\begin{align} 
\alpha_j(\tau,z) 
 = 
 e^{-i m_j\tau}
 \left( 1+\frac{1}{i m_j}\frac{\pa}{\pa \tau} \right)v_j(\tau,z),
 \label{5_3}
\end{align}
$\overline{\alpha_j}$ denotes the complex conjugate of $\alpha_j$, 
and 
$\gamma_j$ is a function of $(\tau,z,\vomega)$ satisfying 
$$
 |\gamma_j(\tau,z,\vomega)|
 \le 
 \frac{C \eps \langle \vomega \rangle^2 e^{4|z|}}{ \tau}, 
$$
for $(\tau, z, \vomega) \in [\tau_0,T)\times \R^2\times \mathbb{H}$. 
In the above estimates, the constants $C$ are independent of 
$\eps$, $T$, $\tau$, $z$, $\vomega$.
\end{lem}

Proof: 
Because of the relations 
$v_k =\Real(\alpha_k e^{im_k\tau})$, 
$\omega_a \pa_{\tau}v_k =-\omega_a m_k \Imag(\alpha_k e^{im_k\tau})$ 
and $m_2=2m_1$, 
we may regard $Q_j(v,\vomega \pa_{\tau}v)$ as a 
trigonometric polynomial in $e^{im_1\tau}$ (with coefficients depending 
on $\alpha_k$, $m_k$, $\vomega$), that is, 
\begin{align}
 Q_j(v,\vomega \pa_{\tau}v)
 = 
 \sum_{\stackrel{1\le k_1 \le k_2 \le 2}{\sigma_1, \sigma_2 \in \{+,-\}}} 
 \Psi_{j k_1 k_2}^{\sigma_1 \sigma_2}(\vomega) 
 \alpha_{k_1}^{(\sigma_1)} \alpha_{k_2}^{(\sigma_2)} 
 e^{i(\sigma_1 k_1+ \sigma_2 k_2) m_1\tau}
\label{5_4}
\end{align}
for $j=1, 2$, where 
$\alpha_k^{(+)}=\alpha_k$, $\alpha_k^{(-)}=\overline{\alpha_k}$ and 
\begin{align*}
 &\Psi_{j k_1 k_2}^{\sigma_1 \sigma_2}(\vomega) =
 \frac{m_1}{2\pi}\int_{0}^{2\pi/m_1} 
 Q_j(\tilde{V}(\theta),\tilde{W}(\vomega,\theta) )
 e^{-i(\sigma_1 k_1+ \sigma_2 k_2) m_1\theta} d\theta
\end{align*}
with $\tilde{V}(\theta) = \bigl(\cos km_1 \theta \bigr)_{k=1,2}$, 
$\tilde{W}(\vomega,\theta) = \bigl( -\omega_a m_k \sin k m_1 \theta 
\bigr)_{{\stackrel{\scriptstyle k=1,2}{\scriptstyle a=0,1,2}}}$. 
Now we focus on the relation $j=\sigma_1 k_1 + \sigma_2 k_2$, 
which implies creation of $e^{im_j\tau}$ in the right hand side of 
(\ref{5_4}). 
We see that this relation is satisfied precisely when 
$(j,k_1,k_2,\sigma_1,\sigma_2)=(1,1,2,-,+)$ or $(2,1,1,+,+)$, 
and that 
$$
 \Psi_{j k_1 k_2}^{\sigma_1 \sigma_2}(\vomega)
 =\left\{
 \begin{array}{cl}
 \Phi_1(\vomega)& \mbox{if } (j,k_1,k_2,\sigma_1,\sigma_2)=(1,1,2,-,+),\\
 \Phi_2(\vomega)& \mbox{if } (j,k_1,k_2,\sigma_1,\sigma_2)=(2,1,1,+,+).\\
 \end{array}
 \right.
$$  
This observation shows that 
$$
 \frac{e^{-im_1 \tau} }{\tau}Q_1(v,\vomega \pa_{\tau}v)
 -
 \frac{\Phi_1(\vomega)}{\tau} \overline{\alpha_1} \alpha_2
$$
and 
$$
 \frac{e^{-im_2 \tau} }{\tau}Q_2(v,\vomega \pa_{\tau}v)
 -
 \frac{\Phi_2(\vomega)}{\tau} \alpha_1^2
$$
are written as sums of the terms in the form
$$
C(\vomega) \alpha_{k_1}^{(\sigma_1)} \alpha_{k_2}^{(\sigma_2)} 
\frac{e^{i\mu \tau}}{\tau} 
$$
with some $\mu \in \R \backslash \{0\}$ and 
$C(\vomega)=O(\langle \vomega \rangle^2)$ ($|\vomega| \to \infty$). 
Eventually we arrive at (\ref{5_1}) and (\ref{5_2}) 
through the identity 
\begin{align*} 
    \alpha_{k_1}^{(\sigma_1)} \alpha_{k_2}^{(\sigma_2)}
    \frac{e^{i\mu \tau}}{\tau}
 =
   \frac{\pa}{\pa \tau} 
   \left( 
   \frac{\alpha_{k_1}^{(\sigma_1)} \alpha_{k_2}^{(\sigma_2)}
   e^{i\mu \tau}}{i\mu\tau}
   \right)
  -
   \frac{\pa}{\pa \tau} 
  \left(  
    \frac{\alpha_{k_1}^{(\sigma_1)} \alpha_{k_2}^{(\sigma_2)}}{\tau} 
  \right)
  \frac{e^{i\mu \tau}}{i\mu}
\end{align*}
combined with the estimates
$$
 |\alpha_k^{(\sigma)}| 
 = 
 \left|v_k \right| + \frac{1}{m_k}\left|\pa_{\tau} v_k \right|
 \le  C \eps^{1/2} e^{2|z|}
$$
and
$$
 |\pa_{\tau} \alpha_k^{(\sigma)}| 
 = 
 \frac{1}{ m_k}
 \bigl| \bigl( \pa_\tau^2 +m_k^2 \bigr)v_k \bigr|
 \le 
  \frac{C \eps^{1/2} e^{2|z|}}{\tau}.
$$
\qed

\section{A lemma on ODE} 
In this section we investigate the behavior as $\tau \gg \tau_0$ of the 
solution $(\beta_1(\tau,z),\beta_2(\tau,z))$ of 
\begin{align}
 \left\{
  \begin{array}{l}
   \displaystyle{i \frac{\pa \beta_1}{\pa \tau} 
   = 
  \frac{\chi_1(z) \Phi_1(\vomega(z))}{\tau} \overline{\beta_1} \beta_2 
  +r_1(\tau,z)}, \\[2mm]
  \displaystyle{i \frac{\pa \beta_2}{\pa \tau}
   = 
  \frac{\chi_2(z) \Phi_2(\vomega(z))}{\tau} {\beta_1}^2
  +r_2(\tau,z)},
  \end{array}
 \right. 
 \quad \tau>\tau_0,
 \label{6_1}
\end{align}
with the initial condition 
\begin{align}
 \sup_{z \in \R^2}\Bigl( |\beta_1(\tau_0,z)|+|\beta_2(\tau_0,z)| \Bigr)
 \le C\eps.
 \label{6_2}
\end{align}
Here 
$\vomega(z)=(\cosh|z|, -z_1\frac{\sinh|z|}{|z|}, -z_2\frac{\sinh|z|}{|z|})$, 
$\Phi_j$ is given by $(\ref{2_1})$, $\chi_j$ is a real-valued function 
satisfying 
$$
 c\le \frac{\chi_1(z)}{\chi_2(z)} \le C
$$ 
with some $C\ge c>0$, and $r_j(\tau,z)$ satisfies 
$$
 \sup_{z\in \R^2}|r_j(\tau, z)| \le \frac{C\eps}{\tau^{2-\delta}}
$$
with some $0<\delta<1$. 
Note that the condition (a) reduces the system $(\ref{6_1})$ to a trivial one, 
that is to say $i\pa_{\tau} \beta_j =O(\eps \tau^{-2+\delta})$, 
so it is easy to see that $(\beta_1,\beta_2)$ stays bounded when $\tau$ 
becomes large. In the following, we will see that a bit weaker assertion is 
valid under the condition (b). 

\begin{lem}\label{lem6_1}
Suppose that the condition (b) is satisfied. 
Let $(\beta_1,\beta_2)$ be the solution of (\ref{6_1})--(\ref{6_2}) 
on $[\tau_0,T)$.Then we have
$$
 \sup_{(\tau, z)\in [\tau_0,T)\times \R^2}
 e^{-2|z|}
 \Bigl(|\beta_1(\tau,z)| + |\beta_2(\tau,z)| \Bigr)
 \le C\eps,
$$
where $C$ is independent of $\eps$, $T$.
\end{lem}

Proof: We first note that both $\Phi_1(\vomega)$ and 
$\Phi_2(\vomega)$ never vanish and that
$$
 |\Phi_1(\vomega) \Phi_2(\vomega)| 
 = \Real\bigl(\Phi_1(\vomega)\Phi_2(\vomega) \bigr)
 = \Phi_1(\vomega)\Phi_2(\vomega)
 \ge C_0
$$
with some strictly positive constant $C_0$ by virtue of (b). We put 
$$
 B_{\eps}(\tau,z)=\Bigl(
  \lambda_1(z)|\beta_1(\tau,z)|^2
  + 
  \lambda_2(z)|\beta_2(\tau,z)|^2 
  +\eps^2
 \Bigr)^{1/2}
$$
with 
$$
 \lambda_1(z)
 =
 e^{-2\langle z \rangle}
 \sqrt{\frac{\chi_2(z)|\Phi_2(\vomega(z))|}{\chi_1(z)|\Phi_1(\vomega(z))|}}, 
 \quad 
 \lambda_2(z)
 =
 e^{-2\langle z \rangle}
 \sqrt{\frac{\chi_1(z)|\Phi_1(\vomega(z))|}{\chi_2(z)|\Phi_2(\vomega(z))|}}.
$$
Then we see that 
\begin{align*} 
\lambda_1(z)
  =
  e^{-2\langle z \rangle}\sqrt{\frac{\chi_2(z)}{\chi_1(z)}} 
 \frac{|\Phi_2(\vomega(z))|}
 {\sqrt{\Real \bigl( \Phi_1(\vomega)\Phi_2(\vomega)\bigr)}}
 \le
 C e^{-2 |z|} (1+|\vomega(z)|^2)
 \le 
 C
\end{align*}
and
$$
 \lambda_2(z) 
 =\frac{e^{-4\langle z \rangle}}{\lambda_1(z)}
 \ge C e^{-4|z|}.
$$
In the same way, we have $\lambda_2(z) \le C$ and 
$\lambda_1(z) \ge C e^{-4|z|}$. 
Therefore 
\begin{align}
 B_{\eps}(\tau,z) 
 \ge 
 C e^{-2|z|} \Bigl(|\beta_1(\tau,z)|+|\beta_2(\tau,z)| \Bigr).
 \label{6_3}
\end{align}
Next we observe that
$$
 \lambda_1(z) \chi_1(z) \Phi_1(\vomega(z)) 
 = 
 \lambda_2(z) \chi_2(z) \overline{\Phi_2(\vomega(z))},
$$
which implies the matrix 
$$
\begin{pmatrix} 
 \lambda_1(z) & 0\\ 
 0 & \lambda_2(z) 
 \end{pmatrix}
\begin{pmatrix} 
 0 & \chi_1(z) \Phi_1(\vomega(z))\beta_1\\ 
 \chi_2(z) \Phi_2(\vomega(z))\overline{\beta_1} & 0 
 \end{pmatrix}
$$
is hermitian. Thus, by rewriting (\ref{6_1}) in the form 
$$
 i\pa_{\tau}\begin{pmatrix} \beta_1\\ \beta_2\end{pmatrix}
 =
 \frac{1}{\tau}
 \begin{pmatrix} 
 0 & \chi_1(z)\Phi_1(\vomega(z))\beta_1\\ 
 \chi_1(z)\Phi_2(\vomega(z))\overline{\beta_1} & 0 
 \end{pmatrix}
 \begin{pmatrix} \beta_1\\ \beta_2\end{pmatrix}
 +
 \begin{pmatrix} r_1\\ r_2\end{pmatrix},
$$
we see that 
\begin{align*}
 B_{\eps}(\tau,z)\pa_{\tau}B_{\eps}(\tau,z)
 &=
 \frac{1}{2}\pa_{\tau}\Bigl(
 \lambda_1 |\beta_1|^2 + \lambda_2 |\beta_2|^2 
 \Bigr)\\
 &=
 \Imag \left\{
 \overline{\begin{pmatrix} \beta_1 & \beta_2\end{pmatrix}} 
 \begin{pmatrix} 
 \lambda_1 & 0\\ 
 0 & \lambda_2 
 \end{pmatrix}
  i\pa_{\tau}  \begin{pmatrix} \beta_1\\ \beta_2\end{pmatrix}
 \right\}\\
  &=
  \Imag \left\{
 \overline{\begin{pmatrix} \beta_1 & \beta_2\end{pmatrix}} 
 \begin{pmatrix} 
 \lambda_1 & 0\\ 
 0 & \lambda_2 
 \end{pmatrix}
 \begin{pmatrix} r_1\\ r_2\end{pmatrix} 
 \right\}\\
 &\le 
 \frac{C\eps}{\tau^{2-\delta}} B_{\eps}(\tau,z).
\end{align*}
Therefore we have 
\begin{align*}
 B_{\eps}(\tau,z) 
 \le 
 B_{\eps}(\tau_0,z) 
 + \int_{\tau_0}^{\infty} \frac{C\eps}{\eta^{2-\delta}} d\eta
 \le
 C\eps,
\end{align*}
which, together with (\ref{6_3}), leads to the desired estimate.
\qed

\section{A priori estimate} 
Now we are in a position to obtain an a priori estimate for the solution of 
(\ref{3_4}), which is the main step of the proof of Theorem \ref{thm2_1}. 
We set 
$$
  M(T):=
 \sup_{(\tau,z)\in [\tau_0,T)\times \R^2}
 e^{-2|z|} \Bigl(
   |v(\tau,z)| + |\pa_{\tau} v(\tau,z)| 
             + \frac{1}{\tau}|\pa_z v(\tau,z)|
 \Bigr) 
$$
for the smooth solution $v=(v_1,v_2)$ to (\ref{3_4}) on 
$\tau \in [\tau_0,T)$. We will prove the following: 
\begin{lem}\label{lem7_1}
There exist $\eps_1>0$ and $C_1>0$ such that 
$M(T)\le \eps^{1/2}$ implies $M(T)\le C_1\eps$ 
for any $\eps \in (0,\eps_1]$. Here $C_1$ is independent of $T$.
\end{lem}

Once this lemma is proved, we can derive global existence of the solution 
in the following way: 
By taking $\eps_0\in (0,\eps_1]$ so that $2C_1 \eps_0^{1/2}\le 1$, 
we deduce that  $M(T)\le \eps^{1/2}$ implies 
$M(T)\le \eps^{1/2}/2$ for any $\eps \in (0,\eps_0]$. Then, by the continuity 
argument, we have $M(T)\le C_1 \eps$ as long as the solution exists. 
Therefore the local solution to (\ref{3_4}) can be extended to the global one. 
Going back to the original variables, we deduce the small data global 
existence for (\ref{1_1})--(\ref{1_2}).

The rest part of this section is devoted to the proof of Lemma \ref{lem7_1}. 
The proof will be divided into two steps: We first derive an auxiliary 
estimate for the energy 
$$
 E_s(\tau) := E_s(\tau;v_1,m_1)+E_s(\tau;v_2,m_2)
$$
under the assumption that $M(T) \le \eps^{1/2}$. Remark that we do not need 
the special structure of the nonlinearity at this stage. Next we will 
prove the improved estimate for $M(T)$ by using the condition (a) or (b). 

\subsection{Energy estimate with moderate growth}

Our goal here is to show $E_{s_1}(\tau) \le C \eps^2 \tau^{\delta}$ 
under the assumption that $M(T) \le \eps^{1/2}$, 
where $s_1 \ge 4$ and $0<\delta<1$. 
We will argue along the same line as \cite{su4}, \cite{su5}. 
We apply (\ref{4_2}) with $s=s_0+s_1+1$ at first, 
where $s_0$ is an integer greater than $2\kappa$. Since Lemma \ref{lem4_3} 
yields 
$$
 \Bigl\|
 \tilde{Q}_j(\tau, \cdot, v, \pa v)\Bigr\|_{(s)}
 \leq 
  \frac{C}{\tau} M(T) E_s(\tau)^{1/2}
 \leq 
  \frac{C}{\tau} \eps^{1/2} E_s(\tau)^{1/2},
$$
we have
$$
  \frac{d}{d \tau} E_{s_0+s_1+1}(\tau)
 \leq 
  \left( 
   \frac{2\kappa + C \eps^{1/2}}{\tau} + \frac{C }{\tau^2} 
  \right) 
  E_{s_0+s_1+1}(\tau)
 \leq 
  \left( 
   \frac{s_0+\frac{1}{2}}{\tau} + \frac{C }{\tau^2} 
  \right) 
  E_{s_0+s_1+1}(\tau).
$$
It follows from the Gronwall lemma that
$$
 E_{s_0+s_1+1}(\tau)
 \leq 
  E_{s_0+s_1+1}(\tau_0) \exp\left(
   \int_{\tau_0}^{\tau} 
   \frac{s_0+\frac{1}{2}}{\eta}+\frac{C}{\eta^2}d\eta 
  \right) 
 \leq 
  C \eps^2 \tau^{s_0+\frac{1}{2}}.
$$
Next, we apply (\ref{4_3}) with $s=s_0+s_1$. Then we have
\begin{align*}
 \frac{d}{d \tau} E_{s_0+s_1}(\tau)
 \le 
   \frac{C}{\tau^{2}}E_{s_0+s_1+1}(\tau)
 + \frac{C \eps}{\tau} E_{s_0+s_1}(\tau)
 \le
   C \eps^2 \tau^{s_0-\frac{3}{2}} + \frac{C\eps}{\tau} E_{s_0+s_1}(\tau),
\end{align*}
which yields
$$
 E_{s_0+s_1}(\tau) \leq C\eps^2 \tau^{s_0-\frac{1}{2} }.
$$
Repeating this procedure recursively, we have
$$
 E_{s_0+s_1+1-n}(\tau) \leq C\eps^2 \tau^{s_0-n+\frac{1}{2} }
$$
for $n=1,2,\cdots,s_0$. Eventually we see that 
$$
 E_{s_1+1}(\tau) \leq C\eps^2 \tau^{1/2}.
$$
Finally, we again use (\ref{4_3}) with $s=s_1$ to obtain
\begin{align*}
 \frac{d}{d \tau} E_{s_1}(\tau)
 \leq 
   \frac{C}{\tau^{2}}E_{s_1+1}(\tau)+ \frac{C \eps}{\tau} E_{s_1}(\tau)
 \leq
   \frac{C \eps^2}{\tau^{3/2}} + \frac{C\eps}{\tau} E_{s_1}(\tau),
\end{align*}
whence we deduce 
$$
 E_{s_1}(\tau)\leq C\eps^2 \tau^{C\eps}
$$ 
for $\tau \in [\tau_0,T)$. 
By choosing $\eps$ so small that $C\eps \leq \delta$, 
we arrive at the desired estimate.

\subsection{Pointwise estimate}
We are going to prove $M(T)\le C\eps$. First we note that
\begin{align*}
 |(\pa_{\tau}^2+m_j^2)v_j|
 =
 \left| \tilde{Q}_j + \frac{1}{\tau^2}\Lambda v_j \right|
 \le 
 \frac{Ce^{(6-\kappa)|z|}}{\tau} M(T)^2 + \frac{C}{\tau^{2}}E_4(\tau)^{1/2}.
\end{align*}
This implies the assumption of Lemma \ref{lem5_1} is satisfied 
if we take $\kappa \ge 6$ and $s_1 \ge 4$. Next we introduce 
$\alpha_j(\tau,z)$ by (\ref{5_3}). Then we see that $\alpha_1$ satisfies 
\begin{align*}
 i \pa_{\tau} \alpha_1
 = 
 \frac{e^{-im_1\tau}}{m_1}(\pa_{\tau}^2+m_1^2)v_1
 =
 \frac{\chi_1(z)\Phi_1(\vomega(z))}{\tau}\overline{\alpha_1}\alpha_2
  + R_1 + \pa_{\tau} S_1,
\end{align*}
where $\chi_1(z)=\chi(z)/m_1$, 
$S_1(\tau,z)=\chi_1(z) \gamma_1(\tau,z,\vomega(z))$ 
with $\gamma_1$ given by Lemma \ref{lem5_1} and
\begin{align*}
 R_1(\tau,z)
 =
  \frac{e^{-im_1\tau}}{m_1}\tilde{Q}_1 -
 \frac{\chi_1(z)\Phi_1(\vomega(z))}{\tau}\overline{\alpha_1}\alpha_2
 - \pa_{\tau} S_1 + \frac{e^{-im_1\tau}}{m_1\tau^2}\Lambda v_1. 
\end{align*}
Since $\tilde{Q}_1$ is given by (\ref{3_3}), it follows from (\ref{5_1}) that 
\begin{align*} 
|R_1(\tau,z)| 
\le& 
 \frac{e^{-\kappa |z|}}{m_1} \left|
 \frac{e^{-i m_1 \tau} }{\tau}Q_1(v,\vomega(z) \pa_{\tau}v)
  -
 \Biggl( 
  \frac{\Phi_1(\vomega(z))}{\tau} \overline{\alpha_1} \alpha_2
  +
  \pa_{\tau}\gamma_1
 \Biggr) \right|\\
 &
 + \frac{Ce^{(6-\kappa)|z|}}{\tau^{2}} M(T)^2
 + \frac{C}{\tau^2} E_4(\tau)^{1/2} \\
\le &
\frac{C\eps \langle \vomega(z)\rangle^2 e^{(4-\kappa)|z|}}{\tau^2} 
+
\frac{C\eps e^{(6-\kappa)|z|}}{\tau^2} 
+
\frac{C\eps}{\tau^{2-\delta}} 
\\
\le &
\frac{C\eps}{\tau^{2-\delta}} 
\end{align*} 
and 
\begin{align*} 
 |S_1(\tau,z)| 
 \le 
 \frac{C\eps e^{(6-\kappa)|z|}}{\tau}
 \le 
 \frac{C\eps}{\tau}.
\end{align*}
Similarly we have
\begin{align*}
 i \pa_{\tau} \alpha_2 
 =
 \frac{\chi_2(z)\Phi_2(\vomega(z))}{\tau}\alpha_1^2
 + R_2 + \pa_{\tau}S_2
\end{align*}
with $\chi_2(z)=\chi(z)/m_2$ and suitable $R_2$, $S_2$ satisfying 
$$
 |R_2(\tau,z)| \le \frac{C\eps}{\tau^{2-\delta}},
 \quad 
 |S_2(\tau,z)| \le \frac{C\eps}{\tau}.
$$
Now, we set $\beta_j=\alpha_j+iS_j$ so that 
$(\beta_1,\beta_2)$ satisfies (\ref{6_1}) with 
\begin{align*} 
 r_1
 =&
  R_1-
 \frac{\chi_1(z)\Phi_1(\vomega(z))}{\tau}
 \bigl(
 i\overline{\alpha}_1 S_2-i\alpha_2\overline{S}_1 + \overline{S}_1S_2
 \bigr),\\
 r_2=&
 R_2 - 
 \frac{\chi_2(z)\Phi_2(\vomega(z))}{\tau} 
 \bigl( 2i\alpha_1 S_1- S_1^2 \bigr).
\end{align*}
Since 
\begin{align*} 
 |r_1|
 \le &
  |R_1|+
  C\frac{e^{(2-\kappa)|z|}}{\tau}
  \bigl(
  |\overline{\alpha}_1|| S_2|
  +
  |\alpha_2||\overline{S}_1| 
  + 
  |\overline{S}_1||S_2|
  \bigr)
  \le 
  \frac{C\eps}{\tau^{2-\delta}},\\
 |r_2|
 \le &
  |R_2|+
  C\frac{e^{(2-\kappa)|z|}}{\tau}
  \bigl( |\alpha_1|| S_1|+|S_1|^2 \bigr)
  \le 
  \frac{C\eps}{\tau^{2-\delta}},
\end{align*}
we can apply Lemma \ref{lem6_1} to obtain 
\begin{align*}
 \sup_{(\tau,z)\in [\tau_0, T)\times \R^2} 
 e^{-2|z|} \bigl| \beta_j(\tau, z) \bigr|
 \le 
 C\eps.
\end{align*}
We thus deduce 
\begin{align*}
  \bigl| \alpha_j(\tau, z) \bigr|
 \le 
  \bigl| \beta_j(\tau, z) \bigr| + \bigl|S_j(\tau,z) \bigr|
 \le 
  C\eps e^{2|z|}.
\end{align*}
Finally, from
$$
  |v_j(\tau, z)| + |\pa_{\tau} v_j(\tau, z)|
  \le
  C \bigl(|v_j(\tau, z)|^2 + m_j^2|\pa_{\tau} v_j(\tau, z)|^2 \bigr)^{1/2}
  = 
   C|\alpha_j(\tau, z)| 
$$
and 
$$
  |\pa_z  v_j(\tau, z) | \le C E_3(\tau)^{1/2}
$$
it follows that 
\begin{align*}
 M(T) 
 \le  
 C\sup_{(\tau,z) \in [\tau_0,T)\times \R^2}e^{-2|z|}\Bigl(
 |\alpha(\tau, z)| +\frac{E_3(\tau)^{1/2}}{\tau}
 \Bigr)
 \le C\eps,
\end{align*}
as desired.\qed

\section{End of the proof of Theorem \ref{thm2_1}} 
The remaining task is to show the decay estimate (\ref{1_3}). 
Remember that our change of variables is 
$$
 u_j(t,x)=\frac{\chi(z)}{\tau} v_j(\tau,z)
$$
with 
$$
 t+2K=\tau \cosh |z|, \quad 
 x_j=\tau \frac{ z_j}{|z|} \sinh |z|
$$ 
for $|x|< t+2K$, and that $u(t,\cdot)$ is supported on 
$\{x \in \R^2 : |x|\le t+K\}$. 
Moreover, we already know that 
$|v(\tau,z)| + |\pa_{\tau} v(\tau,z)| \le C \eps e^{2|z|}$ 
and $|\pa_z v(\tau,z)| \le C \eps \tau^{\delta/2}$. 
So it follows that 
\begin{align*}
 |u(t,x)| 
 =
 \frac{\chi(z)\cosh|z|}{\tau \cosh|z|}  |v(\tau,z)|
 \le
 \frac{e^{-(\kappa -3)|z|}}{t+2K} \cdot e^{-2|z|}|v(\tau,z)|
 \le 
 \frac{C\eps}{1+t}.
\end{align*}
Also, by using (\ref{3_1}), we see that 
\begin{align*}
 \pa_a u(t,x)
 =&
 \frac{\omega_a(z) \chi(z)}{\tau} \pa_{\tau} v(\tau,z)
 +
 \frac{1}{\tau}\sum_{j=1}^{2} \chi(z)\eta_{aj}(z) 
 \frac{\pa_{z_j}v(\tau,z)}{\tau}\\
 & 
 +
 \frac{1}{\tau^2} 
 \Biggl\{
 - \chi(z) \omega_a(z) 
 + \sum_{j=1}^{2} \bigl(\pa_{z_j}\chi(z) \bigr) \eta_{aj}(z)
 \Biggr\}
 v(\tau,z),
\end{align*}
whence
\begin{align*}
 \sum_{|I|=1}|\pa_{t,x}^{I} u(t,x)|
 \le
 \frac{C\eps e^{-(\kappa-3)|z|}}{t+2K} 
 +
 \frac{C \eps e^{-(\kappa-5)|z|}}{(t+2K)^2} 
 \le 
 \frac{C\eps}{1+t}.
\end{align*}
To sum up, we obtain (\ref{1_3}) with $p=\infty$. 
As for the case of $p \in [2,\infty)$, we have
\begin{align*}
  \sum_{|I|\le 1}\|\pa_{t,x}^{I} u(t,\cdot)\|_{L^p(\R^2)}
  =&
 \sum_{|I|\le 1}\|\pa_{t,x}^{I} u(t,\cdot)\|_{L^p(\{x\in \R^2; |x|\le t+K\})}
 \\
 \le&
 C \sum_{|I|\le 1}\|\pa_{t,x}^{I} u(t,\cdot)\|_{L^{\infty}} \cdot 
 \left(\int_{\{x\in \R^2: |x|<t+K\}} \!\! 1 \ \ dx \right)^{1/p}
 \\
 \le &
 C\eps (1+t)^{-1+2/p}, 
\end{align*}
which completes the proof. \qed

\subsection*{Acknowledgement.}
One of the authors (H.S.) is partially supported by Grant-in-Aid for 
Young Scientists (B) (No.18740066 and No. 22740089), MEXT.



\begin{thebibliography}{99}

\bibitem{de} J.-M.Delort,  
{\em 
 Existence globale et comportement asymptotique pour l'\'equation de 
 Klein-Gordon quasi lin\'eaire \`a donn\'ees petites en dimension $1$,
} 
Ann. Sci. \'Ecole Norm. Sup. (4) {\bf 34} (2001), 1--61; 
Erratum, ibid. {\bf 39} (2006), 335--345.


\bibitem{dfx} J.-M. Delort, D. Fang and R. Xue, 
{\em 
 Global existence of small solutions for quadratic quasilinear 
 Klein-Gordon systems in two space dimensions,
} 
J. Funct. Anal. {\bf 211} (2004), 288--323.




\bibitem{hnb} N.Hayashi, P.I.Naumkin and Ratno Bagus Edy Wibowo, 
{\em 
Nonlinear scattering for a system of nonlinear Klein-Gordon equations,} 
J. Math. Phys. {\bf 49} (2008), 103501, 24pp.



\bibitem{kos} S. Katayama, T. Ozawa and H. Sunagawa, 
{\em 
 A note on the null condition for quadratic nonlinear
Klein-Gordon systems in two space dimensions,
} 
preprint, 2011.



\bibitem{su1} H. Sunagawa,
{\em 
 On global small amplitude solutions to systems of cubic nonlinear 
 Klein-Gordon equations with different mass terms in one space dimension,
} 
J. Differential Equations 
{\bf 192} (2003), 308--325.


\bibitem{su2} H. Sunagawa,
{\em 
 A note on the large time asymptotics for a system of Klein-Gordon 
 equations,
} 
Hokkaido Math. J. {\bf 33} (2004), 457--472.


\bibitem{su3} H. Sunagawa, 
{\em
Large time asymptotics  of solutions to nonlinear Klein-Gordon systems,
} 
Osaka J. Math. {\bf 42} (2005), 65--83.


\bibitem{su4} H. Sunagawa,
{\em 
Remarks on the asymptotic behavior of the cubic nonlinear Klein-Gordon 
equations in one space dimension, 
} 
Differential Integral Equations {\bf 18} (2005), 481--494.


\bibitem{su5} H. Sunagawa,
{\em 
 Large time behavior of solutions to the Klein-Gordon equation with 
 nonlinear dissipative terms,
} 
J. Math. Soc. Japan, {\bf 58} (2006), 379--400.


\bibitem{taf} E. Taflin,
{\em 
Simple non-linear Klein-Gordon equations in two space dimensions, with 
long-range scattering, 
}
Lett. Math. Phys.  {\bf 79} (2007), 175--192.


\bibitem{tsu} Y. Tsutsumi,
{\em 
 Stability of constant equilibrium for the Maxwell-Higgs equations, 
} 
Funkcial. Ekvac. {\bf 46} (2003), 41--62.


\end{thebibliography}
\end{document}